\documentclass[conference,a4paper]{IEEEtran}
\usepackage{color}

\usepackage[version=4]{mhchem} % This is for writing CO2 as \ce{CO2}

\usepackage{balance} % This is to allign columns in the last page

\usepackage{empheq} % This is for arrays of equations using "empheq"

\usepackage{cite}
\usepackage[hidelinks]{hyperref}

  \usepackage{graphicx}
  \DeclareGraphicsExtensions{.eps}

\begin{document}

% Copyright:
% \makeatletter
% %%%%%%%%%%%%%%%%%%%%%%%%%%%%%% User specified LaTeX commands.
% \def\ps@IEEEtitlepagestyle{%
%   \def\@oddfoot{\mycopyrightnotice}%
%   \def\@evenfoot{}%
% }
% \def\mycopyrightnotice{%
%   {\footnotesize 978-1-5386-1953-7/17/\$31.00 \textcopyright2017 IEEE\hfill}% <--- Change here
%   \gdef\mycopyrightnotice{}% just in case
% }

\title{Optimal Scheduling of Frequency Services Considering a Variable Largest-Power-Infeed-Loss}
%\title{A Unit Commitment Framework Considering a Variable Largest Power Infeed Loss}
%\title{Flexibility from Nuclear Plants to\\Reduce the Cost of Ancillary Services}
%\textcolor{red}{RUN A Word SPELL CHECK FOR THE LATEX PDF FILE, IT'S BETTER THAN OVERLEAF}

\author{\IEEEauthorblockN{Luis Badesa, Fei Teng and Goran Strbac}
\IEEEauthorblockA{Department of Electrical and Electronic Engineering\\
Imperial College London\\
London SW7 2AZ, United Kingdom\\
Email: \{luis.badesa, f.teng, g.strbac\}@imperial.ac.uk}}

% make the title area
\maketitle

% As a general rule, do not put math, special symbols or citations
% in the abstract
\begin{abstract}
% Limit to 150 words for conferences. From 150 to 200 words for journals.
% It has to include:
% 	- The primary objective of the paper
% 	- Describe your research design and methodology (the methods and procedures you employed) 
% 	- The main outcomes and results
% 	- The conclusions that might be drawn from these data and results. 	
% 	- Include any implications for further research or application/practice.

Low levels of inertia due to increasing renewable penetration bring several challenges, such as the higher need for Primary Frequency Response (PFR). %In turn, this increases the cost for provision of ancillary services. 
A potential solution to mitigate this problem consists on reducing the largest possible power loss in the grid. This paper develops a novel modelling framework to analyse the benefits of such approach.

A new frequency-constrained Stochastic Unit Commitment (SUC) is proposed here, which allows to dynamically reduce the largest possible loss in the optimisation problem. Furthermore, the effect of load damping is included by means of an approximation, while its effect is typically neglected in previous frequency-secured-UC studies. Through several case studies, we demonstrate that reducing the largest loss could significantly decrease operational cost and carbon emissions in the future Great Britain's grid.
\end{abstract}


\section*{Nomenclature}
\addcontentsline{toc}{section}{Nomenclature}

\subsection*{Indices and Sets}
\begin{IEEEdescription}[\IEEEusemathlabelsep\IEEEsetlabelwidth{$\textrm{RoCoF}_{\textrm{max}}$}]
\item[$g,\,\, \mathcal{G}$] Index, Set of thermal generators.
%\item[$\mathcal{M}$] Set of must-run generators.
%\item[$n,\,\, \mathcal{N}$] Index, Set of nodes.
\end{IEEEdescription}

\subsection*{Constants}
\begin{IEEEdescription}[\IEEEusemathlabelsep\IEEEsetlabelwidth{$\textrm{RoCoF}_{\textrm{max}}$}]
%\item[$\Delta\tau(n)$] Time-step corresponding to node $n$ (h).
\item[$\Delta f_{\textrm{max}}$] Maximum admissible frequency deviation from nominal value (Hz).
\item[$\Delta f^{\textrm{ss}}_{\textrm{max}}$] Maximum admissible frequency deviation at quasi-steady-state (Hz).
%\item[$\pi(n)$] Probability of reaching node $n$.
%\item[$\textrm{c}^{\textrm{m}}_g$] Marginal cost of thermal unit $g$ (\pounds/MWh).
%\item[$\textrm{c}^{\textrm{nl}}_g$] No-load cost of thermal unit $g$ (\pounds/h).
%\item[$\textrm{c}^{\textrm{st}}_g$] Startup cost of thermal unit $g$ (\pounds).
\item[$\textrm{D}$] Load damping rate (1/Hz).
\item[$f_0$] Nominal frequency (Hz).
\item[$\textrm{H}_g$] Inertia constant of thermal unit $g$ (s).
\vspace{3pt}
\item[$\textrm{H}^\textrm{L}$] Inertia constant of the generator producing power $P^\textrm{L}$ (s).
\item[$\textrm{P}^{\textrm{D}}$] Total demand (MW).
\item[$\textrm{P}_g^{\textrm{max}}$] Maximum generation of thermal unit $g$ (MW).
\vspace{-7pt}
\item[$\textrm{P}_{\textrm{max}}^\textrm{L}$] Maximum generation of largest unit (MW).
\vspace{4pt}
\item[$\textrm{P}^\textrm{L}_\textrm{i}$] Segment $i$ in the discretisation of $P^\textrm{L}$ in the nadir constraint (MW).
\item[$\textrm{RoCoF}_{\textrm{max}}$] Maximum RoCoF admissible (Hz/s).
\item[$\textrm{T}_\textrm{d}$] Delivery time of PFR (s).
\end{IEEEdescription}

% \subsection*{Semi-Constants \normalfont{(fixed with respect to linear program but variable between time-steps in the UC)}}
% \begin{IEEEdescription}[\IEEEusemathlabelsep\IEEEsetlabelwidth{$\textrm{RoCoF}_{\textrm{max}}$}]

% \end{IEEEdescription}

\subsection*{Decision Variables}
\begin{IEEEdescription}[\IEEEusemathlabelsep\IEEEsetlabelwidth{$\textrm{RoCoF}_{\textrm{max}}$}]
\item[$m^\textrm{L}_i$] Binary variables for discretisation of $P^\textrm{L}$ in the nadir constraint.
\item[$P_g$] Power produced by generator $g$ (MW).
\item[$P^\textrm{L}$] Largest possible power loss (MW).
\item[$P^\textrm{L}_\textrm{nadir}$] Auxiliary variable for the discretisation of $P^\textrm{L}$ in the nadir constraint (MW).
\item[$R_g$] PFR provision from generator $g$ (MW).
\item[$x_g$] Binary variable corresponding to the on/off state of generator $g$.
%\item[$R_g(n)$] Primary Frequency Response provision from thermal unit $g$ at node $n$ (MW).
\end{IEEEdescription}

\subsection*{Linear Expressions \normalfont{(linear combinations of decision variables)}}
\begin{IEEEdescription}[\IEEEusemathlabelsep\IEEEsetlabelwidth{$\textrm{RoCoF}_{\textrm{max}}$}]
%\item[$C_g(n)$] Operating cost of thermal unit $g$ at node $n$ (\pounds).
\item[$H$] System's inertia after generation loss $P_\textrm{L}$ (MW$\cdot \textrm{s}^2$).
%\item[$N_g^{sg}(n)$] Number of thermal unit $g$ that start generating at node $n$.
%\item[$N_g^{up}(n)$] Number of thermal unit $g$ that are online at node $n$.
\item[$R$] Total PFR from all generators (MW).
\end{IEEEdescription}

% For peer review papers, you can put extra information on the cover
% page as needed:
% \ifCLASSOPTIONpeerreview
% \begin{center} \bfseries EDICS Category: 3-BBND \end{center}
% \fi
%
% For peerreview papers, this IEEEtran command inserts a page break and
% creates the second title. It will be ignored for other modes.
\IEEEpeerreviewmaketitle

\section{Introduction}

Integration of Renewable Energy Systems (RES) poses significant problems for grid operators. One of the biggest challenges is due to the reduced level of inertia caused by renewables, which compromises the frequency security of the system. Inertia and Frequency Response (FR) are two services which allow to contain electric-frequency excursions after a power outage. Therefore, the low level of inertia in decarbonised grids greatly increases the need for FR, as demonstrated in \cite{FeiISGT2017} for Great Britain's (GB) system. In turn, this higher need for FR increases both the operational cost of the system and RES curtailment.

One potential solution for this challenge consists on reducing the largest possible outage in the system, as proposed by National Grid in a recent report \cite{NationalGridPLossInertia}. In GB's system, this would be achieved by deloading nuclear plants under certain system's conditions, as these are the largest sources of power in the grid. However, the effectiveness of this option has not yet been analysed, due to the lack of a tool allowing to carry out this study. The present paper focuses on developing such a tool and using it to analyse the operational benefits that this ``deloading approach" would bring, both in terms of reduction in cost and in carbon emissions.

The tool developed here is a Unit Commitment (UC) which optimally schedules inertia and FR, while considering a variable largest-power-infeed-loss. Since inertia and FR are mainly provided by thermal generators, they are inherently related to energy production, and therefore must be scheduled by a UC algorithm. Previous studies such as \cite{UCRestrepoGaliana,LinearizedUC,UCFaroe} have focused on constraining the UC in order to optimally provide inertia and FR. However, to the best of our knowledge, the size of the largest contingency has not yet been modelled in a frequency-secured UC formulation. Some related work was carried out in \cite{OMalleyDeload}, which considered a variable largest outage in a competitive market-dispatch framework. Simulations were used to deduce the frequency-security constraints, an approach which only allows to cover some of the operational conditions. Instead, here we deduce analytical constraints from a mathematical model of the time-evolution of frequency. 

The present paper builds on the work in \cite{FeiStochastic}, which introduced inertia-dependent FR requirements in an SUC. The SUC formulation in \cite{FeiStochastic} is expanded here in order to explicitly model a variable largest-possible-power-loss in the system. The largest loss is now considered as a decision variable in the UC, which then optimises the system's operation by dynamically balancing the cost of deloading and the savings from a reduced need for inertia and FR. Decreasing the largest loss might be optimal depending on the demand and RES generation in the system. 

In addition, here we propose an approximation for the effect of load damping on frequency nadir. Most of the previous work on frequency-secured UC has neglected the effect of load damping, as it would yield complicated mathematical expressions for the nadir requirement. Authors typically argue that the impact of load damping on the need for inertia and FR would be small in any case, and therefore it is acceptable to ignore it. However, \cite{FeiStochastic} demonstrated that a damping factor of 1\%/Hz would reduce operational costs by 5\%, for the 2030 GB power system. The approximation proposed here allows to model the non-negligible effect of load damping, while still giving a simple linear expression for the nadir constraint.

The structure of the paper is as follows: the proposed frequency-constrained UC model is described in Section \ref{SectionUC}. Section \ref{SectionCaseStudies} presents the results of several case studies, demonstrating the value of reducing the largest possible loss in GB's system. Finally, Section \ref{SectionConclusion} gives the conclusion.

\section{UC with Frequency Security Constraints} \label{SectionUC}

%\subsection{Stochastic Unit Commitment}

The UC problem is solved in this paper by using an expanded version of the stochastic scheduling model described in \cite{FeiStochastic}. This SUC is formulated as a Mixed-Integer Linear Program (MILP) which minimises the expected operational cost of the system, while taking into account the uncertainty introduced by wind power. As compared to \cite{FeiStochastic}, where the largest possible loss took a fixed value, the largest loss is modelled here as a decision variable defined as:
\begin{equation} \label{PLossDVdefinition}
P^{\textrm{L}} \geq P_g \qquad \forall g \in \mathcal{G}
\end{equation}
% \textcolor{red}{Be careful with this. Is $P_g=0$ for units that are off? I think it doesn't work because $P_g$ doesn't necessarily have to be 0 for off units, the algorithm can give it any value since it is multiplied by the binary on/off variable, which will be 0 in this case.}
% 
% I could then use this formulation:
% \begin{equation}
% P_{\textrm{L}} = P_g\cdot x_g \qquad  \forall g
% \end{equation}
% % Where x_g is the ON/OFF status of unit g, which is a binary variable and therefore the product can be linearised using big-M
% But I decided not to put this in the paper so that I don't make things more confusing by adding yet another linearisation.
%
% % If I want to consider a model in which just nuclear plants are the largest in the system:
% \begin{equation}
% P^{\textrm{L}} \geq P_g \qquad \forall g \in \mathcal{M}
% \end{equation}
%
% % If units were clustered, it would be:
% \begin{equation}
% P^{\textrm{L}} \geq \frac{P_\mathcal{M}}{\textrm{N}_\mathcal{M}} 
% \end{equation}
% Valid for a model in which equal generators are clustered and nuclear plants are always online.
%
% \textcolor{red}{ACTUALLY, DOES THIS WORK FOR ALL types of units, not only must-run? Isn't $P_g=0$ for units that are off? I think it doesn't work because $P_g$ doesn't necessarily have to be 0 for off units, the algorithm can give it any value since it is multiplied by the binary on/off variable, which is the DV that will be 0 when the unit is off.}
Note that (\ref{PLossDVdefinition}) is easily generalizable to consider any source of power production, such as an interconnector importing power from another grid.

This SUC model assures frequency security by optimally scheduling inertia and FR. %As these two services are mainly provided by the thermal plants producing energy at a particular time, the UC must be constrained so that a certain level of inertia and FR is maintained. 
The deduction of the frequency-security constraints, as well as some linearisations needed for their implementation in an MILP formulation, are given in the following subsections.
%The level of inertia and FR which assures frequency security depends on the wind-demand conditions, \textcolor{red}{as shown in the next subsections}.

% \begin{figure}[!t]
% \centering
% \includegraphics[width=3.5in]{ScenarioTree}
% \caption{Example of a scenario tree used in the Stochastic Unit Commitment}
% \label{ScenarioTree}
% \end{figure}

\subsection{Frequency-Security Constraints} \label{SectionDynamicFrequency}

Frequency security is assured if three requirements are respected: 1) RoCoF must be lower than a certain limit at all times; 2) the frequency nadir must not be below a predefined threshold; and 3) frequency must recover to a certain value 60 seconds after an outage (called ``quasi-steady-state requirement") \cite{NationalGridRequirements}. Certain constraints must be enforced in the UC so that these frequency requirements are met. As explained in \cite{FeiStochastic}, these frequency-security constraints can be deduced from the swing equation, which describes the time evolution of frequency deviation after a generation outage \cite{KundurBook}:
\begin{equation} \label{SwingEq}
2H\frac{\textrm{d} \Delta f(t)}{\textrm{d} t}+\textrm{D}\cdot \textrm{P}^{\textrm{D}}\cdot\Delta f(t)=\sum_{g \in \mathcal{G}}\Delta P_g(t)-P^{\textrm{L}}
\end{equation}
% I don't include subindex "s" in the sum, as Fei does in paper "Stochastic...", because I am not considering any storage in this paper
%Eq. (\ref{SwingEq}) considers the uniform frequency model: all generators are considered to move coherently as a single lumped mass and load damping is represented by a single constant. 
PFR provision by thermal unit $g$ is modelled as:
\begin{equation} \label{PFRdefinition}
\Delta P_g(t)=\left\{ 
\begin{array}{ll}
\frac{R_g}{\textrm{T}_{\textrm{d}}}\cdot t \quad & \mbox{if $t<\textrm{T}_{\textrm{d}}$} \\
R_g \quad & \mbox{if $t\geq \textrm{T}_{\textrm{d}}$}
\end{array}
\right.
\end{equation}
Term $\Delta P_g(t)$ only considers PFR, rather than including secondary and tertiary FR. In this paper, the focus is put on PFR, as the need for PFR is most affected by a low level of inertia. Note that, for the sake of simplicity, PFR is assumed in (\ref{PFRdefinition}) to start being provided right after the generation outage, i.e., the frequency deadband of turbine governors is neglected. However, the frequency deadband could easily be included in the model presented here.

By solving the swing equation (\ref{SwingEq}), the RoCoF, nadir and q-s-s constraints can be obtained (refer to \cite{FeiStochastic} for a detailed description of the mathematical process):
\begin{equation} \label{RocofConstraint}
H \geq \left| \frac{P^\textrm{L}}{2\cdot \mbox{RoCoF}_{\textrm{max}}} \right|
\end{equation}
\begin{multline} \label{FeiKlog}
\frac{2 H\cdot R}{\textrm{T}_{\textrm{d}}}\log\left(\frac{2 H\cdot R}{\textrm{T}_{\textrm{d}}\cdot P^{\textrm{L}}\cdot \textrm{D}\cdot \textrm{P}^{\textrm{D}}+2H\cdot R}\right) \\
\geq (\textrm{D}\cdot \textrm{P}^{\textrm{D}})^2\cdot \Delta f_{\textrm{max}}- P^{\textrm{L}}\cdot \textrm{D}\cdot \textrm{P}_{\textrm{D}}
\end{multline}
\begin{equation} \label{qssConstraint}
R \geq P^{\textrm{L}} - \textrm{D}\cdot \textrm{P}^{\textrm{D}} \cdot \Delta f^{\textrm{ss}}_{\textrm{max}}
\end{equation}
In (\ref{RocofConstraint}), the system's level of inertia after the largest possible outage is given by:
\begin{equation}
H=\frac{\sum_{g \in \mathcal{G}}\textrm{H}_g\cdot \textrm{P}_g^{\textrm{max}}\cdot x_g-\textrm{P}_{\textrm{max}}^\textrm{L}\cdot \textrm{H}^\textrm{L}}{f_0}
\end{equation}

Constraints (\ref{RocofConstraint}) and (\ref{qssConstraint}) are linear and therefore can be directly implemented in an MILP. However, as $P^{\textrm{L}}$ is modelled as a decision variable in the present work, the nadir constraint (\ref{FeiKlog}) becomes nonconvex. It is not possible to linearise this constraint for its inclusion in an MILP, given the logarithmic function and the several bilinear terms $H\cdot R$, some of which appear within the argument of this logarithmic function. Therefore, here we deduce a new nadir constraint, which is obtained again by solving the swing equation (\ref{SwingEq}), but this time neglecting its load-damping term. The nadir constraint then becomes:
\begin{equation} \label{nadirConstraintSimple}
H\cdot R \geq \frac{(P^{\textrm{L}})^2 \cdot \textrm{T}_{\textrm{d}}}{4 \cdot \Delta f_{\textrm{max}}} 
\end{equation}

Neglecting the effect of load damping yields a more conservative nadir constraint, as load damping helps complying with the nadir requirement. However, in the following subsection we propose a linear approximation to the effect of load damping on supporting the frequency nadir. In addition, we linearise the nadir constraint for its inclusion in an MILP.

\subsection{Impact of Load Damping and Linearisation of the Nadir Constraint}

In order to model the effect of load damping on the nadir constraint, we propose the following linear term:
\begin{equation} \label{nadirConstraint}
H\cdot R \geq \frac{(P^{\textrm{L}})^2 \cdot \textrm{T}_{\textrm{d}}}{4 \cdot \Delta f_{\textrm{max}}} - \frac{\textrm{D} \cdot \textrm{P}_{\textrm{D}} \cdot \textrm{T}_{\textrm{d}}}{4} \cdot P^{\textrm{L}} 
\end{equation}

% Note: eq. (\ref{kFei}) has been obtained assuming that $\Delta f_{DB} \approx 0$.

% \begin{figure}[!t]
% \centering
% \includegraphics[width=2.3in]{test2}
% \caption{Feasible regions defined by constraint (\ref{FeiKlog}), for two different values of $\textrm{P}^L$. The feasible region for each value of $\textrm{P}^L$ is the epigraph of each curve.}
% \label{FigNadirLinear}
% \end{figure}

\begin{figure}[!t]
\centering
\includegraphics[width=2.3in]{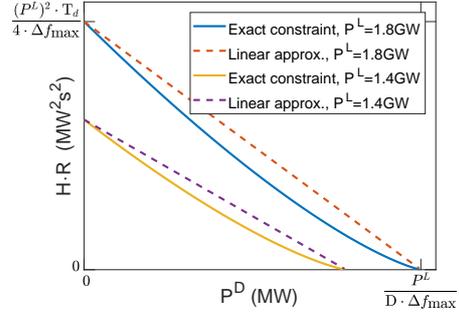}
\caption{Feasible regions defined by constraint (\ref{FeiKlog}), for two different values of $P^\textrm{L}$. The feasible region for each value of $P^\textrm{L}$ is the epigraph of each curve.}
\label{FigNadirLinear}
\end{figure}

This linear term can be deduced by careful examination of (\ref{FeiKlog}). The graphical solution for the exact nadir constraint (\ref{FeiKlog}) is given in Fig. \ref{FigNadirLinear}, for two different fixed values of $P^{\textrm{L}}$. As $H\cdot R=f(\textrm{P}^{\textrm{D}})$ is a convex and monotonically decreasing function, it can be inner-approximated by a line. The y-intercept of that line is given by the right-hand side of (\ref{nadirConstraintSimple}), while the slope of the line can be obtained by considering the largest possible value that $\textrm{P}^{\textrm{D}}$ can take. % As $H\cdot R$ and $\textrm{P}_{\textrm{D}}$ are both positive physical magnitudes, the range of $\textrm{P}_{\textrm{D}}$ for which eq. (\ref{FeiKlog}) has a solution is $\textrm{P}_{\textrm{D}} \in (0,\frac{P_{\textrm{L}}}{\textrm{D}\cdot \Delta f_{\textrm{max}}})$.
Therefore, the linearised effect of load damping on nadir, represented by the dashed lines in Fig. \ref{FigNadirLinear}, is given by (\ref{nadirConstraint}).

Note that the inner approximation of the nadir constraint by a line implies an underestimation of the actual effect of load damping, as can be clearly observed in Fig. \ref{FigNadirLinear}. However, this underestimation is still less conservative than simply neglecting the effect of damping. 

Constraint (\ref{nadirConstraint}) must be linearised before being implemented in an MILP, as it contains two nonlinear terms: $H \cdot R$ and $(P^{\textrm{L}})^2$. In an SUC problem, it is critical to use an MILP formulation: a Mixed-Integer NonLinear Program would considerably increase the computational time needed to solve the problem, so the SUC would likely become intractable. In order to linearise the squared term in constraint (\ref{nadirConstraint}), decision variable $P^{\textrm{L}}$ can be discretised in ``s" segments as follows:
\begin{subequations} \label{NadirArray}
\begin{empheq}[left = \empheqlbrace]{align}
  & H\cdot R +\beta \cdot  \textrm{P}_1^{\textrm{L}} \nonumber \\
  & \qquad \qquad \geq \left(1-m_1^{\textrm{L}} \dotsb - m_{\textrm{s}}^{\textrm{L}} \right) \frac{(\textrm{P}^{\textrm{L}}_1)^2 \cdot \textrm{T}_{\textrm{d}}}{4 \cdot \Delta f_{\textrm{max}}} \tag{\ref{NadirArray}.1} \\
  & H\cdot R +\beta\cdot \textrm{P}_2^{\textrm{L}} \nonumber \\
  & \qquad \qquad \geq \left(1-m_2^{\textrm{L}} \dotsb - m_{\textrm{s}}^{\textrm{L}} \right) \frac{(\textrm{P}^{\textrm{L}}_2)^2 \cdot \textrm{T}_{\textrm{d}}}{4 \cdot \Delta f_{\textrm{max}}} \tag{\ref{NadirArray}.2} \\
  & \qquad \mathrel{\makebox[\widthof{=}]{\vdots}} \nonumber \\
  & H\cdot R +\beta\cdot \textrm{P}_{\textrm{s-1}}^{\textrm{L}} \geq \left(1-m_{\textrm{s}}^{\textrm{L}}\right) \frac{(\textrm{P}^{\textrm{L}}_{\textrm{s-1}})^2 \cdot \textrm{T}_{\textrm{d}}}{4 \cdot \Delta f_{\textrm{max}}} \tag{\ref{NadirArray}.(s-1)} \\
  & H\cdot R +\beta \cdot \textrm{P}_{\textrm{s}}^{\textrm{L}} \geq \frac{(\textrm{P}^{\textrm{L}}_{\textrm{s}})^2 \cdot \textrm{T}_{\textrm{d}}}{4 \cdot \Delta f_{\textrm{max}}} \tag{\ref{NadirArray}.s} \\
  & m_1^{\textrm{L}}+m_2^{\textrm{L}}+ \dotsb + m_{\textrm{s}}^{\textrm{L}} \leq 1 \tag{\ref{NadirArray}.(s+1)}
\end{empheq}
\end{subequations}
Where $\beta$ is defined as:
\begin{equation}
\beta = \frac{\textrm{D}\cdot \textrm{P}^{\textrm{D}} \cdot \textrm{T}_{\textrm{d}}}{4}
\end{equation}

The binary variables $m_i^{\textrm{L}}$ enforce that only one of the constraints in (\ref{NadirArray}) is activated. For this discretisation of $P^{\textrm{L}}$ to hold, a new decision variable must be defined:
\begin{equation} \label{PLossNadir}
P_{\textrm{nadir}}^{\textrm{L}} = m_1^{\textrm{L}} \cdot \textrm{P}_1^{\textrm{L}} + m_2^{\textrm{L}} \cdot \textrm{P}_2^{\textrm{L}} \dotsb + m_{\textrm{s}}^{\textrm{L}} \cdot \textrm{P}^{\textrm{L}}_{\textrm{s}}
\end{equation}
The following constraint applies to $P_{\textrm{nadir}}^{\textrm{L}}$:
\begin{equation} \label{PLossNadirConstraint}
P_{\textrm{nadir}}^{\textrm{L}} \geq P^{\textrm{L}}
\end{equation}

By enforcing constraints (\ref{NadirArray}) and (\ref{PLossNadirConstraint}), the squared term $(P^{\textrm{L}})^2$ in (\ref{nadirConstraint}) is linearised. Then, the bilinear term $H\cdot R$ appearing in each of the constraints in (\ref{NadirArray}) must also be linearised, for which we use a big-M method as proposed in \cite{FeiStochastic}. With these two linearisations, the nadir constraint can be implemented in an MILP.

%\textcolor{red}{NOTE: I chose to say that $H\cdot R$ is a bilinear term that can be linearised using a big-M, instead of saying the we actually use the linearisation in Fei's paper ``Full SUC''. That's because I never mention in this paper that units are clustered, so it's much simpler to explain the linearisation using big-M. About not mentioning that units are clustered, it doesn't make sense to mention it here, since I don't mention either the computational time that it takes ot run the simulations, and clustering is just to reduce computational time.}

\section{Case Studies} \label{SectionCaseStudies}
In order to demonstrate the benefits from dynamically limiting the largest power infeed loss, the SUC model described in section \ref{SectionUC} was used to run several case studies. Each case study simulated one full year of operation of the 2030 GB power grid. The results of these simulations are analysed here, in terms of the operational cost of the system, load factor of large nuclear units and \ce{CO2} emissions.

The characteristics of the 2030 GB system used as the platform for our simulations can be found in Table I of \cite{LuisISGT2017}. The load damping factor, $\textrm{D}$, was set to 1\%/Hz. For the SUC, a scenario tree branching only in the current-time node was used, and net-demand quantiles of 0.005, 0.1, 0.3, 0.5, 0.7, 0.9 and 0.995 were considered (refer to \cite{AlexEfficient} for further explanation on scenario trees in SUC). %A penalty for emissions of 150 \pounds /ton$\mbox{CO}_{2}$ was also considered.
Different wind-penetration levels were analysed, as the amount of renewable generation that will be present in the GB system by 2030 is still uncertain. As shown in coming subsections, reducing the largest possible loss has a higher value for increasing wind penetration.

%The characteristics of this system are given in Table \ref{Table2030GB}. 

% \begin{table}[!t]
% % increase table row spacing, adjust to taste
% \renewcommand{\arraystretch}{1.2}
% \caption{Characteristics of Thermal Plants}
% \label{Table2030GB}
% \centering
% % Some packages, such as MDW tools, offer better commands for making tables
% % than the plain LaTeX2e tabular which is used here.
% \begin{tabular}{l| l l l}
% \hline
% & Nuclear & CCGT & OCGT\\
% \hline
% Number of Units & 6 & 110 & 30\\
% Rated Power (MW) & 1800 & 500 & 200\\
% Min Stable Gen (MW) & 1800 & 200 & 50\\
% No-Load Cost (\pounds/h) & 0 & 7809 & 8000\\
% Marginal Cost (\pounds/MWh) & 10 & 51 & 110\\
% Startup Cost (\pounds/) & n/a & 9000 & 0\\
% Startup Time (h) & n/a & 4 & 0\\
% Min Up Time (h) & n/a & 4 & 0\\
% Min Down Time (h) & n/a & 1 & 0\\
% Inertia Constant (s) & 5 & 5 & 5\\
% Max Response (MW) & 0 & 50 & 20\\
% Response Slope & 0 & 0.5 & 0.5\\
% Emissions (kg$\mbox{CO}_{2}$/MWh) & 0 & 368 & 833\\
% \hline
% \end{tabular}
% \end{table}

As mentioned before, reducing the largest possible power loss could be achieved in GB's system by deloading large nuclear plants. It is important to remark that the reason behind deloading nuclear units might not always be to reduce the cost of frequency services, i.e., the cost of providing inertia and PFR. In high-wind-generation conditions, deloading nuclear allows to accommodate more wind, which is zero-cost energy, therefore reducing the total cost of energy provision. However, the present study focuses on deloading nuclear plants just to reduce the need for inertia and FR.

\subsection{Cost of Providing Frequency Services} \label{CostFreq}

\begin{figure}[!t]
\centering
\includegraphics[width=2.3in]{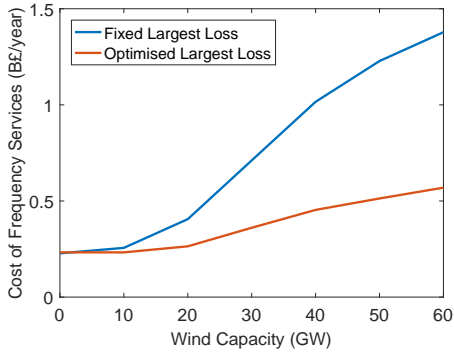}
\caption{Annual cost of frequency services under different wind-penetration scenarios, for a largest possible loss of 1.32GW.}
\label{FigCostFreq}
\end{figure}

\begin{figure}[!t]
\centering
\includegraphics[width=2.3in]{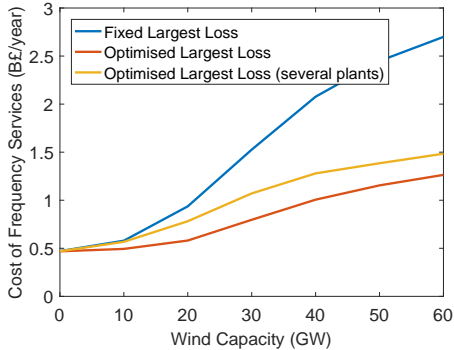}
\caption{Annual cost of frequency services under different wind-penetration scenarios, for a largest possible loss of 1.8GW.}
\label{FigNuclearFleet}
\end{figure}

In this subsection we analyse the cost of providing frequency services, namely inertia and FR. The amount of inertia and FR needed to comply with the frequency-security constraints is provided by running part-loaded thermal generators. Running a high number of part-loaded generators is more expensive than producing the same amount of energy from a lower number of fully-loaded generators, and it potentially causes RES curtailment. The difference in operational cost between these two cases is what we refer to as ``cost of frequency services".

For these simulations, we compare three different scenarios: the current largest loss in GB's system, of 1.32GW; the projected largest loss in 2030, which will be of 1.8GW; and a variation of the latter scenario in which several nuclear plants, not just one, are rated at 1.8GW.

First of all, we analyse the benefits from deloading based on the current single largest plant in GB. In Fig. \ref{FigCostFreq} the annual cost of frequency services for the 1.32GW-largest-loss case is presented. Two different operational strategies are considered: in ``Fixed Largest Loss", the largest nuclear plant is forced to operate at maximum output at all times; in ``Optimised Largest Loss", this plant is allowed to reduce its power output by 33\% of its rating. %In addition, two different RoCoF requirements are studied: a strict requirement of $\textrm{RoCoF}_\textrm{max}=0.125 \textrm{Hz/s}$ and a relaxed requirement of 0.5Hz/s. Relaxing the RoCoF requirement has been proposed by National Grid \cite{NationalGridRelaxRocof}.
Although the quantitative results presented in Fig. \ref{FigCostFreq} might vary depending on the characteristics of the power system studied, this figure shows a clear trend that would hold in any system: reducing the largest loss has a significant positive impact in the cost of frequency services, particularly for high-wind-penetration cases. 

Here we also study the impact of a higher largest-possible-loss in the future GB system. In Fig. \ref{FigNuclearFleet} we present the cost of frequency services for the projected largest loss of 1.8GW, as well as a variation of this scenario in which 6 nuclear plants have a 1.8GW rating. %We consider only the projected relaxed RoCoF requirement ($\textrm{RoCoF}_\textrm{max}=0.5 \textrm{Hz/s}$), and assume that 
Large nuclear plants are again allowed to deload by 33\% of its rating. Note that the cost of frequency services for the ``Fixed Largest Loss" case is the same for both the scenario with one large plant rated at 1.8GW, and for the one with 6 large plants. By comparing this ``Fixed Largest Loss" case in Fig. \ref{FigNuclearFleet} with the one in Fig. \ref{FigCostFreq}, one can notice that the cost of frequency services doubles Fig. \ref{FigNuclearFleet} for every wind-penetration scenario. This issue should be brought to the attention of system planners, since a larger nuclear plant may be a sensible option from the energy-efficiency point of view, but its implications in increased operational cost of the system must also be considered. Fig. \ref{FigNuclearFleet} also demonstrates that deloading large nuclear units brings even further savings, in absolute terms, for a 1.8GW-largest-loss when compared to the 1.32GW case. Regarding the scenario with 6 plants rated at 1.8GW, as all 6 plants must be deloaded in order to effectively reduce the largest possible loss, the deloading strategy is less effective, although it still leads to significant savings.

Finally, one can notice that for all scenarios in both Fig. \ref{FigCostFreq} and \ref{FigNuclearFleet}, the cost of frequency services increases with increasing wind penetration, as would be expected: when non-synchronous wind generation, which does not contribute to inertia or FR, replaces conventional generators, the system's levels of inertia and FR are reduced; then, part-loaded conventional generators must be brought online just to provide frequency services, therefore increasing the operational cost of the system.

\subsection{Analysis of the Load Factor of Nuclear Units}

In subsection \ref{CostFreq} the savings in operational cost due to reducing the largest possible loss have been demonstrated. However, the impact of deloading nuclear units on the investment return of these generation plants should also be considered. 

Nuclear plants have very high investments costs, and they are expected to provide inexpensive, carbon-free, base-load energy. However, if these plants do not operate in a base-load mode, but are deloaded in order to reduce the largest possible loss, the investment might be less attractive. While the present work focuses on the operational aspects of the power grid, and therefore investment costs are not taken into account, we analyse the load factor of the largest nuclear unit to inform system planners of this issue.

\begin{figure}[!t]
\centering
\includegraphics[width=2.4in]{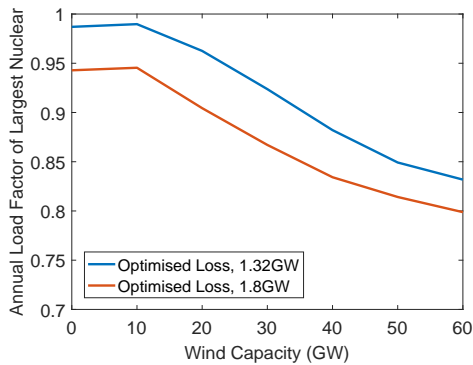}
\caption{Annual load factor of the largest nuclear unit as a function of wind capacity. Two different ratings of the largest nuclear unit are considered.}
\label{FigLoadFactor}
\end{figure}

Fig. \ref{FigLoadFactor} presents the results of our study. One can notice that %RoCoF regulation has a big impact on the load factor of the largest nuclear plant. In addition, 
an increased rating of the largest nuclear plant decreases the load factor, as the 1.8GW plant is more frequently deloaded than the 1.32GW plant. If deloading nuclear units were a strategy to be implemented by system operators, system planners should be aware that this strategy might lead to a significantly reduced load factor for these units.

\balance

\vspace{3mm}
\subsection{Impact on \ce{CO2} Emissions}

\begin{figure}[!t]
\centering
\includegraphics[width=2.3in]{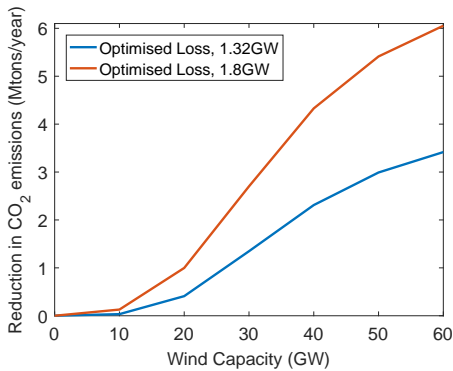}
\caption{Reduction in carbon emissions due to allowing the largest nuclear unit to deload, as a function of wind power capacity.}
\label{FigEmissions}
\end{figure}

As mentioned in subsection \ref{CostFreq}, in certain occasions a number of part-loaded thermal generators must be online in order to provide inertia and FR. In high-wind-generation conditions, wind power might be curtailed in order to keep these part-loaded thermal plants online. By reducing the largest possible loss, less inertia and FR is needed, and therefore less part-loaded plants are required to be online. Then, by reducing the largest loss, carbon emissions are also reduced, since more wind power can be accommodated. 

Fig. \ref{FigEmissions} shows the amount of carbon emissions that would be cut annually by allowing the largest nuclear plant to deload. A very significant reduction in \ce{CO2} emissions is achieved, particularly for cases of high wind capacity and a 1.8GW largest loss. Nevertheless, even for a 1.32GW-largest-loss case the reduction in emissions is considerable. Given the strict emission targets recently set in countries all over the world, reducing the largest possible loss has proved to be an effective strategy to comply with this legislation.

\vspace{10mm}
\section{Conclusion} \label{SectionConclusion}

This paper has introduced a frequency-constrained UC in which a variable largest-possible-power-loss is explicitly modelled. This UC model has been used to analyse the potential operational benefits from deloading large nuclear generators in GB's system. These benefits, which have been demonstrated to be considerable, particularly in a high-wind-penetration scenario, are both economic and in terms of a reduction in carbon emissions.

%As compared to much of the previous work on frequency-secured UC, the present paper does not neglect the effect of load damping. A linear approximation to the effect of load damping on the nadir constraint has been proposed, so that this constraint can be included in an MILP formulation.

The UC framework presented here can be used to support the discussion on different options to tackle the frequency-security problem in the low-inertia system. In the future, this framework could be extended to consider other services such as Enhanced Frequency Response. In addition, further work on studying the operational benefits from deloading nuclear plants should focus on analysing its interaction with other frequency services. It would be particularly interesting to study its interaction with extra inertia provision, a service already considered in \cite{LuisISGT2017}. This is motivated by National Grid's recent report \cite{NationalGridPLossInertia}, where it was stated that reducing the largest loss would be more effective than creating a market for inertia. However, these two strategies could be complimentary, so it should be determined if it is an economically sensible option to make them coexist.

Finally, it would be critical to analyse the implication on investment return of a partially-loaded nuclear plant, as the results presented here have demonstrated that the deloading strategy could significantly reduce the load factor of large nuclear units.

% % use section* for acknowledgment
% \section*{Acknowledgment}

% trigger a \newpage just before the given reference
% number - used to balance the columns on the last page
% adjust value as needed - may need to be readjusted if
% the document is modified later
%\IEEEtriggeratref{8}
% The "triggered" command can be changed if desired:
%\IEEEtriggercmd{\enlargethispage{-5in}}

\bibliographystyle{IEEEtran} 
\bibliography{Luis_PhD}

\end{document}